\documentclass[11pt,a4paper]{article}
\usepackage[T1]{fontenc}
\usepackage{lmodern}
\usepackage{amsmath,amssymb,amsthm}
\usepackage[margin=30mm,top=26mm,bottom=29mm]{geometry}
\usepackage{microtype,tikz,needspace}
\usepackage[hidelinks]{hyperref}
\hypersetup{pdftitle={A uniform additive bound for rainbow girth in Aharoni’s conjecture},
  pdfauthor={Yiming Hao, Hui Lei, Jianing Liu, Zhijun Xu},
  pdfkeywords={rainbow cycle, edge-colored graph, directed girth, Caccetta-Häggkvist conjecture}}
\newtheorem{theorem}{Theorem}[section]
\newtheorem{lemma}[theorem]{Lemma}
\newtheorem{conjecture}{Conjecture}
\newcommand{\rg}{\operatorname{rg}}
\newcommand{\authormail}[1]{\mbox{\urlstyle{same}\href{mailto:#1}{\nolinkurl{#1}}}}
\newcommand{\mathschool}{School of Mathematical Sciences and LPMC, Nankai University, Tianjin 300071, China.}
\newcommand{\frontlabel}[1]{\par\addvspace{.6\baselineskip}\small\noindent
  \hspace*{1.35em}\textbf{#1.}\enspace\ignorespaces}
\newenvironment{keywords}{\frontlabel{Key words}}{\par}

\makeatletter
\renewcommand{\@maketitle}{\begin{center}
  \vspace*{0.5em}
  {\LARGE\bfseries\@title\par}
  \vspace{0.65em}
  {\normalsize\scshape\@author\par}
\end{center}}
\makeatother
\title{A uniform additive bound for rainbow girth in Aharoni’s conjecture}
\author{Yiming Hao\thanks{\mathschool\ Email: \authormail{1120230031@mail.nankai.edu.cn}.},
  Hui Lei\thanks{School of Statistics and Data Science, LPMC and KLMDASR, Nankai University, Tianjin 300071, China. Email: \authormail{hlei@nankai.edu.cn}.},
  Jianing Liu\thanks{\mathschool\ Email: \authormail{liujianing@mail.nankai.edu.cn}.},
  and Zhijun Xu\thanks{\mathschool\ Email: \authormail{zhj_xu@mail.nankai.edu.cn}.}}
\date{}
\begin{document}
\maketitle
\begin{abstract}
The rainbow girth of an edge-colored graph is the minimum length of a cycle whose edges
have pairwise distinct colors. Aharoni's generalization of the Caccetta--H\"aggkvist
conjecture asserts that a simple $n$-vertex graph with $n$ color classes, each of size
at least $r$, has rainbow girth at most $\lceil n/r\rceil$. Hompe and Huynh proved the
upper bound $n/r+\alpha_r$, where $\alpha_r=O(r^5\log^2 r)$. We answer their question about the dependence of the additive
term on $r$ by proving that the rainbow girth is at most $n/r+5000$ for all $n$ and $r$.
\end{abstract}

\section{Introduction}
A simple digraph is loopless and has at most one arc from $u$ to $v$ for each ordered
pair of distinct vertices $u,v$. Throughout, $n$ and $r$ are positive integers. In 1978, Caccetta and H\"aggkvist proposed the following conjecture on short directed cycles.

\begin{conjecture}[Caccetta and H\"aggkvist~\cite{CH}]\label{conj:CH}
Let $D$ be a simple $n$-vertex digraph with minimum out-degree at least $r$.
Then $D$ contains a directed cycle of length at most $\lceil n/r\rceil$.
\end{conjecture}

For each fixed $r$, the conjecture holds for all sufficiently large $n$.
Indeed, Shen~\cite{ShenGirth} proved that it holds whenever $n\ge 2r^2-3r+1$.
A second line of work bounds the error uniformly in $r$. Chv\'atal and
Szemer\'edi~\cite{CS} proved that $D$ contains a directed cycle of length at most
$n/r+2500$. Nishimura~\cite{Nishimura} reduced
the additive constant to $304$, and Shen~\cite{Shen} reduced it to $73$. We use Shen's
additive bound in the proof of our main theorem.

In 2017,  Aharoni proposed an edge-coloring version of Conjecture~\ref{conj:CH},
later recorded by Aharoni, DeVos and Holzman~\cite{ADH}. Let $G$ be a finite simple
undirected graph. A graph is called \emph{simple} if it has no loops or parallel edges. An \emph{edge-coloring} assigns a color to each edge, and the edges
assigned the same color form a \emph{color class}. Adjacent edges may have the same color.
A subgraph is \emph{rainbow} if its edges have pairwise distinct colors. The \emph{rainbow
girth} $\rg(G)$ is the length of a shortest rainbow cycle, with $\rg(G)=\infty$ if no such
cycle exists. 

\begin{conjecture}[Aharoni~\cite{ADH}]\label{conj:aharoni}
If an $n$-vertex simple graph has $n$ color classes, each of size at least $r\ge1$, then
$\rg(G)\le\lceil n/r\rceil$.
\end{conjecture}

To see why Conjecture~\ref{conj:aharoni} implies Conjecture~\ref{conj:CH}, consider a
digraph $D$ as in the latter. A pair of opposite arcs gives a directed cycle of length two,
which suffices since $r\le n-1$. Otherwise, forget the orientations and color each edge by
the tail of its arc. Each vertex is the tail of at least $r$ arcs, so each color class
has at least $r$ edges. On a rainbow
cycle of length $k$, the $k$ tails are distinct vertices of the cycle. Thus each vertex has
one outgoing and one incoming arc on the cycle, so the cycle is directed. Conversely, every
directed cycle is rainbow under this coloring.

The same correspondence gives examples attaining the conjectured bound~\cite{ADH}.
For $n>2r$, take the digraph on $\mathbb Z/n\mathbb Z$ with arcs from $i$ to
$i+1,\ldots,i+r$. The steps on any directed cycle sum to a positive multiple of $n$, so
the cycle has length at least $k=\lceil n/r\rceil$. Taking $k-1$ steps of size $r$ and a
final step of size $n-(k-1)r$ gives a cycle of length $k$. coloring by tails therefore
gives a simple graph with $r$ edges of each color and rainbow girth exactly $k$.

The case $r=1$ of Conjecture~\ref{conj:aharoni} follows from the one-edge-per-color
argument above. DeVos et al.~\cite{DeVos} proved the sharp bound $\lceil n/2\rceil$ for
$r=2$. The exact bound for $r=3$ remains open. In that case Clinch, Goerner, Huynh and
Illingworth~\cite{CGHI} proved that $\rg(G)\le 4n/9+7$. For general $r$, Hompe and
Spirkl~\cite{HS} showed that $\rg(G)\le n/r$ if every color class has at least
$10^{11}r$ edges. Hompe and Huynh~\cite{HH} subsequently obtained an additive approximation to Aharoni's conjecture,
with a correction term for colour classes of size less than $r$.

\begin{theorem}[Hompe and Huynh {\cite[Theorem~1.7]{HH}}]\label{thm:HH-defect}
For every integer $r\ge2$, there exists a constant $\alpha_r=O(r^5\log^2 r)$ such that
the following holds. Let $G$ be a simple graph on $n$ vertices whose edges are partitioned
into $n$ colour classes $E_1,\ldots,E_n$. If $2\le |E_i|\le r$ for each $i$, then
\[
    \rg(G)\le \frac{n+\sum_{i=1}^n (r-|E_i|)}{r}+\alpha_r .
\]
\end{theorem}

If every colour class has at least $r\ge2$ edges, retain exactly $r$ edges of each colour.
Applying Theorem~\ref{thm:HH-defect} to the resulting spanning subgraph gives
$\rg(G)\le n/r+\alpha_r$, since the sum vanishes. The case $r=1$ holds with $\alpha_1=0$,
so this bound holds for all $r\ge1$~\cite[Corollary~1.8]{HH}. Hompe and Huynh asked how
the additive term in this corollary depends on $r$. Our main theorem shows that it can
be chosen independently of $r$.

\begin{theorem}\label{thm:main}
Let $n$ and $r$ be positive integers. Every simple graph on $n$ vertices with $n$ edge colors and
at least $r$ edges of each color contains a rainbow cycle of length at most $n/r+5000$.
\end{theorem}

In the cyclic example every color class is a star, meaning that all its edges have a
common endpoint. Guo~\cite[Theorem~8]{Guo} showed that for each fixed $\varepsilon>0$, if
every class has at least two edges and at least $\varepsilon n$ of the $n$ classes are
non-stars, then $\rg(G)=O_{\varepsilon}(\log n)$. In our proof the edges of monochromatic
stars supply the forest to be contracted.

We henceforth keep exactly $r$ edges of each color and put $x=n/r$. Since
$nr\le\binom{n}{2}$, we have $x>2$. We write $\log$ for the natural logarithm and
$\log_2$ for the logarithm to base two. If we can delete $q$ vertices while keeping an
edge of every color, then choosing one edge of each color gives $n$ edges on $n-q$
vertices. Contracting a forest preserves this excess of $q$ edges. A short cycle after
contraction may, however, pass through a long path in one of the trees when lifted back
to $G$. In Section~\ref{sec:weighted} we bound the length of such a cycle by assigning
weights to the surviving edges according to the distances of their endpoints from the
roots. In Section~\ref{sec:main} we choose the vertices to delete and apply the weighted
bound. When $r$ is large compared with $x$, we either apply Shen's bound directly or
choose the roots first so that every representative edge meets a root. In the latter
case a cycle of length $k$ after contraction lifts to a rainbow cycle of length at most $2k$.

\section{Short cycles after contracting a forest}\label{sec:weighted}

\begin{figure}[ht]
\centering
\begin{tikzpicture}[x=0.92cm,y=0.70cm,every node/.style={font=\small},
    point/.style={circle,fill=black,inner sep=1.55pt}, tree/.style={thin}, oncycle/.style={line width=.9pt}]
  \draw[densely dashed,gray,rounded corners=3pt] (.65,-.57) rectangle (3.35,1.89);
  \node[gray] at (.93,1.47) {$T$};
  \draw[oncycle] (0,0)--(1,0)--(2,0.7)--(3,0)--(4,0);
  \draw[tree] (2,0.7)--(2,1.6);
  \foreach \p in {(0,0),(1,0),(2,0.7),(2,1.6),(3,0),(4,0)} \node[point] at \p {};
  \node[below] at (1,0) {$y$}; \node[below] at (3,0) {$z$};
  \node[right] at (2,1.6) {$s$};
  \draw[->] (4.6,0.45)--(6.0,0.45) node[midway,above] {contract};
  \draw[oncycle] (6.7,0.45)--(9.3,0.45);
  \foreach \p in {(6.7,.45),(8,.45),(9.3,.45)} \node[point] at \p {};
  \node[below] at (8,.45) {$T$};
\end{tikzpicture}
\caption{Part of a cycle through the tree $T$ inside the dashed box. The root is $s$. Its path from $y$ to $z$ uses two tree edges, while $\ell(y)+\ell(z)=4$.}
\label{fig:tree}
\end{figure}
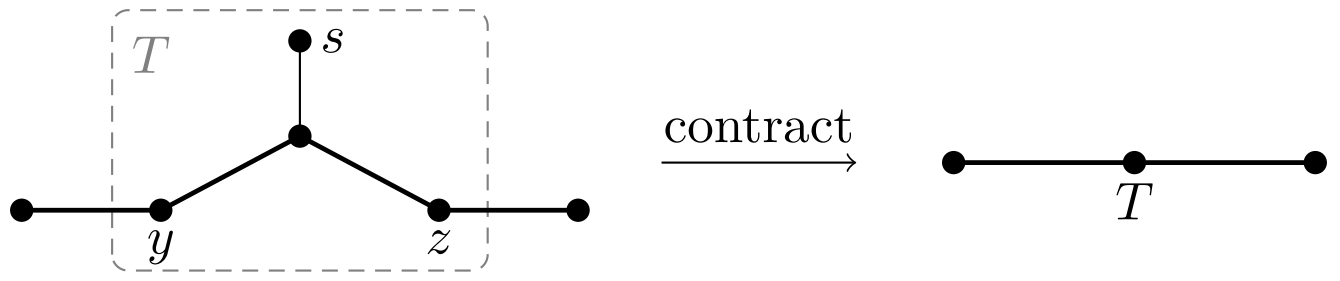

Let $F$ be an acyclic spanning subgraph of a rainbow graph $R$. Choose a root in each tree and
write $\ell(v)$ for the distance from $v$ to its root. Contract each tree to a single vertex. We
keep all edges outside $F$, including any that become loops or parallel edges. Give the image of
each such edge $yz$ weight
\begin{equation}\label{eq:lift}
                    1+\ell(y)+\ell(z).
\end{equation}
A cycle in the contracted graph gives a rainbow cycle in $R$ of length at most its weight. To
recover the cycle, join its two attachment points in each tree by the unique path between them.
The length of this path is at most the sum of their distances to the root. A tree occurs only once
on the cycle, and distinct trees are disjoint. This also proves the assertion for a loop or a pair
of parallel edges, which we count as cycles of length one and two, respectively.
Figure~\ref{fig:tree} illustrates the length bound inside one tree.

We next bound the weight of a shortest cycle. As in the work of Bollob\'as and
Szemer\'edi~\cite{BS}, the relevant parameter is the \emph{excess} of a graph, its number of edges
minus its number of vertices. For $q\ge1$, write
\[
   w(q)=\frac{2\log_2(q+1)}{3q}.
\]
The function $w$ is decreasing, $w(1)=2/3$, and $w(q)\le \log(q+1)/q$.

\begin{lemma}\label{lem:weighted-cycle}
Let $q\ge1$ be an integer. If a multigraph of excess at least $q$ has nonnegative edge weights of
total weight $W$, then it has a cycle of weight at most $w(q)W$.
\end{lemma}

\begin{proof}
First join the components by bridges of weight zero. Keep a spanning tree and $q+1$ edges outside
it, so the remaining graph is connected and has excess $q$. Repeatedly delete leaves and replace
paths whose internal vertices have degree two by single edges of the same total weight. The
resulting graph has minimum degree at least three, since its excess is positive. Replace each
vertex of degree $s>3$ by a path on $s-2$ vertices, give the path edges weight zero, and
distribute the original edge ends along the path so that each new vertex has degree three. These
operations preserve the excess. The resulting cubic multigraph therefore has $2q$ vertices and
$3q$ edges. A cycle in it yields a cycle of at most the same weight in the original graph. Indeed,
contracting the added paths and expanding the suppressed paths gives a nonempty closed trail, and
no bridge lies on such a trail.

It remains to prove the bound for a cubic multigraph. We may assume its weights are positive and
rational, since the minimum of the weights of its finitely many cycles is continuous. Let $W'$ be
its total weight and $g$ its minimum cycle weight. For any sufficiently large multiple $D$ of the
common denominator, replace an edge of weight $a$ by a path of $Da$ unit edges. The subdivided
graph has $DW'$ edges, $2q$ vertices of degree three, and girth $Dg$.

Choose a uniformly random oriented unit edge and start at its midpoint. Walk to its head, choose
an outgoing edge uniformly among those other than the reverse of the incoming edge, and walk to
the next midpoint. Repeat this step $h=\lfloor Dg/2\rfloor$ times. At every step the oriented edge
is uniform. Indeed, each orientation leaving a vertex of degree $s$ has $s-1$ possible
predecessors, each followed by that orientation with probability $1/(s-1)$. There are two choices
at a vertex of degree three and one at a vertex of degree two. If $K$ counts these steps, then
\[
   \mathbb E K = h\,\frac{3(2q)}{2DW'}=\frac{3qh}{DW'}.
\]

Fix the initial midpoint and consider the trees of possible walks in its two directions. View the
initial midpoint as a subdivision vertex. All branching vertices occur at distance at most
$h-\tfrac12<Dg/2$ along the walk. Each such vertex is reached by a unique route from the midpoint.
A repeated vertex on one route would give a cycle, and two distinct simple routes to the same
vertex would give a cycle of length less than $Dg$ in their union. Thus the two trees together
have at most $2q$ branching vertices and at most $2q+2$ leaves. A leaf reached after $k$ choices
has conditional probability $2^{-k-1}$, including the choice of direction. Since its contribution
to the conditional expectation of $2^K$ is $1/2$, $\mathbb E 2^K\le q+1$, also after averaging
over midpoints. By convexity, $2^{\mathbb E K}\le \mathbb E 2^K$. We obtain
$3qh/(DW')\le\log_2(q+1)$. Letting $D$ tend to infinity gives $g\le w(q)W'$, as required.
\end{proof}

We now choose the forest. Choosing more roots shortens the paths in the forest, but it also leaves
more edges outside the forest. We estimate their total weight. In a simple edge-colored graph
$H$, let $U$ be a set of vertices. At each $u\in U$, choose a monochromatic star centred at $u$,
using a different color at each centre. We call these the \emph{designated stars}. Use their
edges, oriented away from their centres, to form a digraph on $V(H)$. Every vertex outside $U$ is
a sink. Each directed cycle is rainbow because its arcs leave distinct centres. For positive
integer vertex weights $\mu(v)$ set $\mu(A)=\sum_{v\in A}\mu(v)$, and let $N^+(u)$ denote the
leaves of the designated star at $u$.

\begin{lemma}\label{lem:star-forest}
Suppose that $H$ has $m$ vertices and $m+q$ nonempty colors, where $q\ge1$. Let $b=m+q-|U|$ be
the number of colors with no designated star. Suppose that $\mu(N^+(u))\ge d$ for every $u\in U$,
where $d\ge1$ is an integer. If $\rg(H)>J$ for an integer $J>146$, set $\theta=1-146/J$. For each
$c>0$, put $M=(1-e^{-c\theta})^{-1}$. Then
\begin{equation}\label{eq:forest-bound}
 \rg(H)\le w(q)\left[
        \frac{cM}{d}\,\mu(U)+|U|(J+1)e^{-c\theta J}+(1+2M)b
                 \right].
\end{equation}
\end{lemma}

\begin{proof}
The three terms in~\eqref{eq:forest-bound} will count the weights of edges from randomly chosen
roots, edges from roots added later, and edges of the other colors. To choose the roots, we first
estimate the total vertex weight reachable from a centre. Let $B_j(u)$ be the set of vertices
reachable from $u$ by a directed path of length at most $j$, including $u$. We claim that for
$1\le j\le J$,
\begin{equation}\label{eq:ball-growth}
 B_j(u)\text{ contains no sink}
 \quad\Longrightarrow\quad \mu(B_j(u)\setminus\{u\})\ge d\theta j .
\end{equation}
Replace every vertex $v$ by $\mu(v)$ copies and each arc by all arcs in the same direction between
the corresponding sets of copies. A directed cycle of length at most $J$ in this digraph would
project to a directed closed walk containing a rainbow cycle of length at most $J$ in $H$.

Start at one copy of $u$ and take the layers at directed distances $1,\dots,j$. None contains a
copy of $u$, since a return to any such copy gives a closed walk in the original digraph. Write
$a$ for the total number of vertices in these layers. In particular
$a\le\mu(B_j(u)\setminus\{u\})$. We use Shen's bound~\cite{Shen} on a digraph made from $k=\lfloor
J/j\rfloor$ copies of the layers arranged cyclically. Within each copy, keep all arcs leaving
layers $1,\dots,j-1$. Their heads lie in the same copy, since one arc increases distance by at
most one. From every vertex in layer $j$, draw arcs to $d$ distinct vertices of the first layer of
the next copy. The first layer has at least $d$ vertices, and the hypothesis that $B_j(u)$
contains no sink gives outdegree at least $d$ at every earlier layer. This digraph is simple. In
particular, if $k=1$ then $j>1$, so the new arcs join different layers.

Any directed cycle using one of the new arcs uses at least $k$ of them. Between two consecutive
new arcs it goes from layer $1$ to layer $j$ and so traverses at least $j-1$ other arcs. A cycle
using only old arcs has length greater than $J$. The directed girth is therefore at least $kj$.
Shen's bound on the $ka$ vertices gives $kj\le ka/d+73$. As $k\ge J/(2j)$, this yields $a\ge
d(j-73/k)\ge d\theta j$, proving~\eqref{eq:ball-growth}.

Choose each centre $v\in U$ independently with probability $1-e^{-c\mu(v)/d}$. Let $S_0$ contain
the chosen centres and all sinks. Add to $S_0$ every vertex whose directed distance to $S_0$ is
greater than $J$, counting an unreachable set as having infinite distance. Call the resulting set
$S$ and let $\ell(v)$ be the directed distance from $v$ to $S$. We have $\ell(v)\le J$. If
$B_j(v)$ contains a sink, it meets $S_0$. Otherwise, by independence and~\eqref{eq:ball-growth},
it misses $S_0$ with probability at most $e^{-c\theta j}$. In particular the probability that a
vertex belongs to $S\setminus S_0$ is at most $e^{-c\theta J}$. Summing the tails gives
\[
    \mathbb E\ell(v)=\sum_{j=0}^{J-1}\Pr(\ell(v)>j)
       \le \sum_{j=0}^{\infty} e^{-c\theta j}=M .
\]

Each vertex outside $S$ chooses an outgoing edge on a shortest path to $S$. These edges form a
forest with one root in $S$ per tree, since the distance decreases by one on each chosen edge. For
each centre $u\in S$ retain the edge $uv$ with $v\in N^+(u)$ of minimum $\ell(v)$. Retain also one
fixed edge of every color not used by a designated star. The forest edges and the retained edges
give exactly one representative of every color. On contracting the forest we obtain a multigraph
with $|S|$ vertices and $|S|+q$ edges. Weight the retained edges by~\eqref{eq:lift}. Here $\ell$
agrees with distance to the root in the forest, so the observation preceding
Lemma~\ref{lem:weighted-cycle} applies.

To bound the total weight, consider first a centre $u$ chosen in the initial sample. Its retained
edge has weight $1+\min_{v\in N^+(u)}\ell(v)$. If $B_j(u)\setminus\{u\}$ meets $S_0$, a path to a
vertex of this intersection starts at $u$ and takes its first step to a leaf. That leaf is within
distance $j-1$ of $S$. The choice of $u$ leaves all other sampling decisions independent. Thus,
conditional on choosing $u$, the probability that the minimum leaf distance is at least $j$ is at
most $e^{-c\theta j}$ for $1\le j\le J$, by~\eqref{eq:ball-growth}. The conditional mean weight is
at most $1+\sum_{j=1}^{J}e^{-c\theta j}\le M$. Since $1-e^{-z}\le z$, the initially chosen centres
contribute at most $cM\mu(U)/d$ in expectation. A centre added in the second step contributes at
most $J+1$. Each of the $b$ other edges has expected weight at most $1+2M$. Adding these bounds
and applying Lemma~\ref{lem:weighted-cycle} to a choice of $S$ whose total weight is at most the
mean proves~\eqref{eq:forest-bound}.
\end{proof}

\section{Proof of the main theorem}\label{sec:main}

There are two ways to keep every color while deleting vertices. For large $r$ we keep a
representative of every color at a small set of roots. In the other range we keep enough leaves
at the centres of monochromatic stars to apply Lemma~\ref{lem:star-forest}.

A set of vertices \emph{covers} a set of edges if it meets every edge in
the set. The following observation will be used for both deletions. It includes the pairs used by
Hompe and Huynh~\cite[Claim~4.1]{HH}.

\needspace{8\baselineskip}
\begin{lemma}\label{lem:pairs}
Let $F$ be a nonempty set of edges and let $N$ be a set of vertices.
If no vertex in $N$ covers $F$ by itself, then there are at most four
pairs of distinct vertices in $N$ such that every subset of $N$
covering $F$ contains one of these pairs.
\end{lemma}

\begin{proof}
Fix an edge $e\in F$. For each endpoint $u\in N$ of $e$, choose an
edge $e_u\in F$ avoiding $u$. Record the pair $\{u,v\}$ for each
endpoint $v\in N$ of $e_u$. This gives at most four pairs. Any cover
contained in $N$ contains an endpoint $u$ of $e$ and an endpoint $v$
of $e_u$, so it contains one of the pairs.
\end{proof}

\begin{lemma}\label{lem:roots}
If $\log r\ge 24+6\log x$, then
$\operatorname{rg}(G)\le\max\{x+75,116\}$.
\end{lemma}

\begin{proof}
Set $a=\lfloor r/(x+2)\rfloor$. Since $x>2$ and $r\ge e^{24}x^6$,
we have $r/(4x)\le a<r/4$. At each vertex incident with at least
$r-a$ edges of one color, choose one such color and designate the
star formed by its edges at that vertex. Distinct centres have distinct
colors, since two such stars of the same color would together contain
at least $2(r-a)-1>r$ edges. Let $N_0$ be the set of vertices without a
designated star and put $t_0=|N_0|$. Exactly $t_0$ colors are undesignated.

If $t_0\le a$, delete $N_0$ and orient the surviving edges of each
designated star away from its centre. Each remaining vertex has
outdegree at least $r-2a$. Every directed cycle is rainbow, so
Shen's theorem~\cite{Shen} gives a cycle of length at most
\[
   \frac{n}{r-2a}+73\le x+75,
\]
where the inequality follows from $a\le r/(x+2)$.

Suppose that $t_0>a$. We shall choose a small set $S$ and delete
vertices outside $S$ so that each color still has an edge incident
with $S$. We require $S$ to meet the leaves of each designated star.
For each undesignated color class $F$ and each $v\in N_0$, we also
require $S$ to meet the set of endpoints of the edges of $F$ avoiding
$v$. Each of these sets
has at least $\sqrt a$ vertices. Indeed, a designated star has at
least $r-a$ leaves. In the second case at least $a$ edges avoid $v$,
and fewer than $\sqrt a$ vertices cannot span $a$ edges in a simple
graph.

There are at most $2n^2$ sets to meet. Sample each vertex independently
with probability $\min\{1,3\log n/\sqrt a\}$, then add one vertex from
each set missed by the sample. Each set is missed with probability at
most $n^{-3}$, so the expected number of additions is at most $2/n$.
Some resulting set $S$ therefore has size
\[
    s=|S|\le 1+3n\log n/\sqrt a \le 4n\log n/\sqrt a.
\]
We shall use the estimate
\begin{equation}\label{eq:rootratio}
   \frac{a}{36s}
      \ge \frac{\sqrt r}{1152x^{5/2}\log n}
      \ge r^{1/24}.
\end{equation}
The first inequality follows from $a\ge r/(4x)$ and $n=rx$.
For the second, $r\ge e^{24}x^6$ gives $x^{5/2}\le e^{-10}r^{5/12}$
and $\log n\le (7/6)\log r\le (28/e)r^{1/24}$. The last bound follows
by maximizing $y e^{-y/24}$ for $y\ge0$, and $28\cdot1152<e^{11}$.

For each undesignated color, restrict its class to the edges incident
with $S$. For every $v\in N_0$ this restricted class has an edge
avoiding $v$, by the choice of $S$. Apply Lemma~\ref{lem:pairs} to
each restricted class, obtaining at most $4t_0$ pairs in total. We can
delete any subset of $N_0\setminus S$ containing none of these pairs.
Each restricted class then retains an edge. Each designated star also
retains an edge to $S$, since its centre lies outside $N_0$.

To find a large such subset, sample each vertex of $N_0\setminus S$
independently with probability $1/8$. For each recorded pair with
both vertices sampled, remove one of them from the sample. The
expected number left is at least $(t_0-s)/8-4t_0/64$. By
\eqref{eq:rootratio} we have $s<t_0/4$, so this expectation exceeds
$t_0/32$. We may thus delete a set of $q\ge t_0/36$ vertices.

Choose one surviving edge of each color incident with $S$, and discard
all other edges. Delete the isolated vertices outside $S$. At each
remaining vertex outside $S$, contract one incident edge. These
contracted edges form a forest of stars rooted in $S$, together with any isolated roots.
The resulting multigraph has $s$ vertices and at least $s+q$ edges.
Indeed, before contraction it has $n$ edges and at most $n-q$ vertices,
and each contraction reduces both counts by one.

Every surviving edge has at least one endpoint in $S$, and the other endpoint is at distance at
most one from its root. Thus~\eqref{eq:lift} assigns weight at most two to every edge of the
contracted graph. A cycle of length $k$ in it gives a rainbow cycle of length at most $2k$ in $G$.

It remains to find a cycle of length at most $58$. We may assume the
multigraph is simple. Repeatedly delete vertices of degree less than
$1+q/s$. The ratio of edges to vertices stays at least $1+q/s$, so a
nonempty subgraph remains, of minimum degree $d\ge1+q/s$. Since
$q\ge t_0/36>a/36$, \eqref{eq:rootratio} gives $d-1\ge r^{1/24}$.
If the girth exceeded $58$, every walk of length $29$ that never
immediately retraces an edge would be a path. Two such paths from
the same vertex to the same endpoint would contain a cycle of length
at most $58$. From any one vertex there would therefore be at least
$(d-1)^{29}\ge r^{29/24}>r^{7/6}\ge n\ge s$ distinct endpoints.
This contradiction proves the lemma.
\end{proof}

We turn to the deletions for the weighted bound. From now on the stars we use are entire color
classes. Assume $r\ge3$, and call a color a \emph{star color} if all its edges
have a common endpoint. At each such endpoint choose one of its star
colors. Let $U$ be the set of these centres, let $N=V(G)\setminus U$,
and put $t=|N|$. Exactly $t$ colors have not been chosen. This count includes any further star
colors at centres where we have already chosen one.

\begin{lemma}\label{lem:delete}
If $t>0$, one can delete a vertex of $N$ and give the remaining
vertices positive integer weights so that every color survives, every
chosen star has leaves of total weight at least $d$, and
\begin{equation}\label{eq:one-deletion}
   \frac{\mu(U)}{d}\le x\left(1+\frac1{\max\{t,r-1\}}\right).
\end{equation}
If $x\le64t$, one can instead delete $q$ vertices of $N$, keep every
color, and leave at least $d$ leaves at every chosen centre, where
\[
\begin{array}{c|cc}
 & q & d\\ \hline
3\le r<64 & \lceil t/48\rceil & \lceil r/4\rceil\\
r\ge64    & \lceil t/17\rceil & \lceil r/2\rceil .
\end{array}
\]
\end{lemma}

\begin{proof}
If a vertex meets every edge of a color, it is the centre of a star color and belongs to $U$. Deleting any one
vertex of $N$ therefore preserves all colors and leaves at least
$r-1$ edges in each chosen star. Unit weights give
$\mu(U)/d\le n/(r-1)=x(1+1/(r-1))$.

For the other bound in \eqref{eq:one-deletion}, delete a uniformly
random vertex $v$ of $N$. Start with unit weights. Whenever a chosen
star loses a leaf, add one to the weight of a surviving leaf,
choosing one in $N$ if possible. This restores total leaf weight at
least $r$. A star can increase $\mu(U)$ only if $v$ was its unique
leaf in $N$, an event of probability at most $1/t$. For some choice
of $v$ we thus have $\mu(U)\le |U|(1+1/t)$. Taking $d=r$ gives the
second bound.

For the larger deletions, sample each vertex of $N$ independently
with probability $p$. We will delete sampled vertices after restoring
enough of them to preserve the colors. For each nonstar color,
take the pairs supplied by Lemma~\ref{lem:pairs}. Restore one endpoint
of each pair whose two endpoints were sampled. There are at most $t$
nonstar colors, so the expected number restored in this way is at
most $4tp^2$.

For every star color, restore enough leaves to keep at least $d$
edges. Its centre lies outside $N$. If $Y$ of its leaves were sampled
and $k=r-d+1$, at most $\binom{Y}{k}$ restorations suffice. Indeed, when
$Y\ge k$, at most $Y-k+1\le\binom{Y}{k}$ leaves are needed.
The expected cost per star is at most $\binom rk p^k$, by counting
$k$-subsets of sampled leaves. All restorations can be chosen from
the initial sample and then combined. As there are at most
$n=rx\le64rt$ stars, the expected number of vertices left to delete
is at least
\begin{equation}\label{eq:restoration}
 t\left(p-4p^2-64r\binom rk p^k\right).
\end{equation}

For $3\le r<64$, use $p=1/32$ and $d=\lceil r/4\rceil$.
The last term in parentheses is at most $1/160$. For $r=3,4,5$ it
is respectively $192/32^3$, $256/32^4$ and $1600/32^4$.
For $r\ge6$, the inequalities $k\ge3r/4$ and $\binom rk\le2^r$
bound it by $64r\,2^{-11r/4}$, which decreases with $r$ and is
less than $1/160$ at $r=6$. Thus \eqref{eq:restoration} exceeds $t/48$.
For $r\ge64$, use $p=1/8$ and $d=\lceil r/2\rceil$. Now $k\ge r/2$,
so the last term is at most $64r\,2^{-r/2}\le4096\,2^{-32}<1/272$.
Since $p-4p^2=1/16$, \eqref{eq:restoration} exceeds $t/17$.
In each range some sample allows the required number of deletions.
\end{proof}

For small $r$ and $t$ we use the estimate
\begin{equation}\label{eq:few-centres}
 \rg(G)\le x+(4r+2)t+2r\qquad(t\ge1).
\end{equation}
In the proof of \cite[Lemma~5.1 and Claim~5.2]{HH} they construct
disjoint sets $T_1,\ldots,T_t$, one for each vertex in $N$. Write $W$
for the vertices outside these sets and $h$ for the number of centres
in $W$ with a chosen edge leaving $W$. Their construction gives numbers
$d_i\ge0$ with $h+\sum_i d_i\le(4r+1)t$. Any two vertices in $T_i$
are joined inside $T_i$ by a rainbow path of chosen colors of length
at most $|T_i|/r+d_i$. Here all color classes have size $r$, so their
terms counting missing edges vanish. If $W$ is empty, contract the
sets $T_i$ and keep one edge of each of the $t$ remaining colors.
The resulting multigraph has $t$ vertices and $t$ edges and hence a
cycle. Joining its edges by the paths inside the sets gives a rainbow
cycle of length at most $t+n/r+\sum_i d_i$. Indeed, two disjoint sets
cannot each contain an edge of the same chosen star. If $W$ is nonempty,
each of its vertices has a chosen edge to another vertex of $W$.
Orient the chosen stars outwards in $W$. Only the $h$ centres can have
outdegree below $r$, and each has outdegree at least one. By
\cite[Corollary~2.4]{HH} there is a directed cycle of length at most
$(|W|+(r-1)h)/r+2r-1$. Either estimate implies \eqref{eq:few-centres}.

\begin{proof}[Proof of Theorem~\ref{thm:main}]
For $r=1$, one edge per color gives $n$ edges on $n$ vertices and
hence a rainbow cycle. For $r=2$, DeVos et al.~\cite{DeVos} give the
bound $\lceil n/2\rceil$. Suppose $r\ge3$. If $t=0$, orient all
chosen stars outwards and apply Shen's theorem~\cite{Shen}. By
Lemma~\ref{lem:roots} we may thus assume $t>0$ and
\begin{equation}\label{eq:parameter-range}
  \log r<24+6\log x.
\end{equation}
If $\max\{t,r-1\}<10$, then $t\le9$ and $r\le10$, and
\eqref{eq:few-centres} already gives $\rg(G)\le x+398$.

Suppose that $\rg(G)>x+5000$ and put $J=\lfloor x+5000\rfloor$.
Then $\theta=1-146/J>0.97$. After deleting $q$ vertices by
Lemma~\ref{lem:delete}, the graph still has $n$ colors, of which
$b=t$ are unchosen. Any rainbow cycle in this graph also lies in $G$,
so we may apply \eqref{eq:forest-bound}. Two estimates will be used
in each application. For $1/4\le c\le1$, the middle term in
\eqref{eq:forest-bound} satisfies
\[
 |U|(J+1)e^{-c\theta J}
 \le x\,2e^{24}J^7 e^{-0.24J}<0.001x.
\]
Indeed, \eqref{eq:parameter-range} gives $r<e^{24}x^6$, while
$x\le J$ and $|U|\le rx$. The coefficient decreases for $J\ge5000$
and is already less than $0.001$ at $5000$.
Also, for $0<z\le1$, the inequality $e^z\ge1+z+z^2/2$ gives
\[
 \frac1{1-e^{-z}}\le 1+\frac1{z+z^2/2}
     =\frac1z+\frac{1+z}{2+z}\le\frac1z+\frac23 .
\]

First delete one vertex. Since $\max\{t,r-1\}\ge10$,
\eqref{eq:one-deletion} gives $\mu(U)/d\le1.1x$.
With $c=1/4$ we have $M<4.8$ and $w(1)=2/3$, so
\[
 \rg(G)\le \tfrac23\bigl[(1.2\cdot1.1+0.001)x+10.6t\bigr]
 < 0.882x+7.1t.
\]
For $t\le700$ the last term $7.1t$ is at most $4970$.
For $x\ge64t$ the whole bound is less than $x$.
It remains to consider $t>700$ and $x<64t$, when the larger
deletions in Lemma~\ref{lem:delete} apply with unit weights.

If $r<64$, delete $q=\lceil t/48\rceil\ge15$ vertices.
Then $\mu(U)/d\le4x$. Take $c=2/5$, so that $M<3.25$. We have
$w(q)\le w(15)=8/45$ and $tw(q)\le48\log(q+1)$. Substitution in
\eqref{eq:forest-bound} gives
\[
 \rg(G)<0.927x+360\log(q+1)\le0.927x+360\log(3x).
\]
Here $q+1\le2+t/48\le2+64x/48\le3x$, since $t\le n=rx$ and
$x>2$. Using $\log(x/5000)\le x/5000-1$, the last bound is at most
$0.999x+360(\log 15000-1)<x+3600$.

If $r\ge64$, delete $q=\lceil t/17\rceil\ge42$ vertices.
Now $\mu(U)/d\le2x$. Take $c=1$, so that $M<1.70$. Since
$43^2<2^{11}$, we have $w(q)\le w(42)<11/126$.
Also $tw(q)\le17\log(q+1)\le17\log n$. It follows from
\eqref{eq:forest-bound} and \eqref{eq:parameter-range} that
\[
  \rg(G)<0.30x+75\log n
       <0.30x+1800+525\log x .
\]
Finally, $\log(x/750)\le x/750-1$ gives
$525\log x\le0.70x+525(\log 750-1)<0.70x+3150$.
Thus $\rg(G)<x+4950$, again a contradiction.
\end{proof}

\section*{Use of AI}
Discussions with ChatGPT helped us explain how to select vertices for
deletion while retaining an edge of every color, and how distances to
the roots of a forest bound the length of a rainbow cycle lifted from
the contracted graph.
We also used ChatGPT to polish the language of the manuscript.
The authors take full responsibility for the mathematical content and
the final manuscript.

\end{document}